\documentclass[a4paper,12pt]{amsart}
\usepackage[margin=3cm]{geometry}
\usepackage{url, amssymb}
\newtheorem{theorem}{Theorem}

\newtheorem{lemma}[theorem]{Lemma}

\newtheorem{corollary}[theorem]{Corollary}
\newcommand{\dquer}{\overline\partial}
\newcommand{\dquers}{\overline\partial ^*_\varphi}

\numberwithin{equation}{section}
\begin{document}
\title{On the resolvent of the Dirac operator in $\Bbb R^2$.}

\author{Klaus Gansberger}

\thanks{Supported by the FWF-grants  P19147 and P19667.}

 \address{K. Gansberger, Institut f\"ur Mathematik, Universit\"at Wien,
Nordbergstrasse 15, A-1090 Wien, Austria}
\email{klaus.gansberger@univie.ac.at}

\maketitle

\begin{abstract} ~\\
In the present paper, we prove an abstract functional analytic criterion for a class of linear partial differential operators acting on a domain $\Omega\subseteq\Bbb R^n$ which are elliptic in the interior to have compact resolvent. This extends known results for magnetic Schr\"{o}dinger operators to more general differential operators. We point out the relationship between the Dirac operator in real dimension two and the $\dquer$-Laplacian on a certain weighted space on $\Bbb C$ and we use this connection to prove a non-compactness result for its resolvent.
\end{abstract}

\section{Introduction and Results.}~\\

The aim of this paper is to prove a non-compacness result for the resolvent of the Dirac operator in $\Bbb R^2$. This essentially shows that if the magnetic field $B$ has a specific sign, i.e. $B(x,y)\ge0 $ or $B(x,y)\le0$, the spectrum of the Dirac operator $\Bbb D$ is never purely discrete.\\ 
In mathematical physics, the Dirac equation models the behavior of a  ``free'' relativistic spin-$\frac{1}{2}$ particle, see e.g. \cite{tha} for an introduction to and details on the physical interpretation. The state space of such an particle is $L^2(\Bbb R^2,\Bbb C^2)$, so $\Bbb D$ acts $a\ priori $ on $\mathcal C_0^\infty(\Bbb R^2)\oplus \mathcal C_0^\infty(\Bbb R^2)$, the space of smooth functions with compact support. Using the standard choice of Pauli matrices $\sigma_j$ 
\begin{align*}
\label{pauli}
\sigma_1=
\begin{pmatrix}
0 &1\\
1 &0\\
\end{pmatrix}
\quad
\text{and}
\quad
\sigma_2=
\begin{pmatrix}
0 &-i\\
i &0\\
\end{pmatrix},
\end{align*}
the Dirac operator in $\Bbb R^2$ is given by
\begin{equation}
\label{dirac}
\Bbb D=\sigma_1\left(-i\frac{\partial}{\partial x}-A_1(x,y)\right)+\sigma_2\left(-i\frac{\partial}{\partial y}-A_2(x,y)\right)
\end{equation}
acting on $\Psi=(\Psi_1,\Psi_2)\in \mathcal C_0^\infty(\Bbb R^2)\oplus \mathcal C_0^\infty(\Bbb R^2)$ by matrix-multiplication. Here, $A_1$ and $A_2$ are multiplication operators by real-valued functions. It is classical, that $\Bbb D$ is essentially self-adjoint and can be extended in a unique way to a self-adjoint operator acting on $L^2(\Bbb R^2,\Bbb C)\oplus L^2(\Bbb R^2,\Bbb C)$, see for instance \cite{tha} or \cite{hnw}.\\
There is also a notion of the Dirac operator in real dimension three, see e.g. \cite{hnw}. In real dimension two, it is conjectured that the Dirac operator never admits a compact resolvent, cf. \cite{er}, \cite{hnw}. If we denote the spectrum of $\Bbb D$ by $\sigma(\Bbb D)$, this means that the operator $(\Bbb D-\lambda)^{-1}$ is not compact for all $\lambda\in\Bbb C\setminus\sigma(\Bbb D)$. The main result of \cite{hnw} is the following: Let the magnetic field be
$$B(x,y)=\frac{\partial A_2}{\partial x} (x,y)-\frac{\partial A_1}{\partial y} (x,y) $$
for smooth functions $A_j, \ j=1,2$ and define
\begin{align*}
m_q (x,y) =\sum_{|\alpha|=q-1}|\partial^\alpha B (x,y) |\qquad\text{and}\qquad m^r (x,y) =1+\sum_{q=1}^rm_q (x,y).
\end{align*}
Suppose that there exists a sequence of pairwise disjoint balls each one of radius greater than $1$, such that 
\begin{equation}
\label{helffer}
m_{r+1} (x,y)\le Cm^r (x,y)
\end{equation}
holds on the union of these balls. Then the Dirac operator has non-compact resolvent.
Note that this condition is for instance satisfied, if the magnetic potentials are polynomials. 
\vskip 0.5 cm

We will use a different approach, actually coming from complex analysis. The first time that a connection between complex analysis and the Dirac operator was noticed, was in \cite{hahe}. Throughout the paper, we will assume that 
$$B (x,y)=\triangle\varphi(x,y)$$ 
for some function $\varphi$, thus $A_1(x,y)=-\varphi_y(x,y)$ and $A_2(x,y)=\varphi_x(x,y)$. This is not as specific as it seems at the first glimpse, since starting with $B$ one can first find a function $\varphi$ such that $\triangle\varphi=B.$ See also \cite{er} and the references therin for this point of view. Nevertheless, we will put the from the physical point of view rather restrictive regularity assumption $\varphi\in\mathcal C^2(\Bbb R^2)$. But we add as a Remark that at least in the case $B(x,y)\ge0$ this can easily be weakened to assuming that there is a (subharmonic) $\mathcal C^2$-function, such that the difference to $\varphi$ is bounded. This can be seen using the arguments of \cite{gh}, in particular Lemma 2.3, combined with the unitary equivalence of the Pauli operators to a complex Laplacian in an $L^2$-space weighted with $e^{-\varphi}$, see \cite{hahe}. The case $B(x,y)\ge0$ corresponds to subharmonicity of $\varphi$, which from the complex analysis point of view is the interesting one.\\
Let us now formulate our main result.

\begin{theorem}
\label{dirac}
Suppose that the magnetic field 
$$B(x,y)=\frac{\partial A_2}{\partial x}(x,y)-\frac{\partial A_1}{\partial y}(x,y)$$
is of the form $B=\triangle\varphi$ for some $\mathcal C^2$-function $\varphi$. Suppose furthermore that $\varphi$ can be chosen to be bounded from above or from below. Then the Dirac operator has non-compact resolvent.
\end{theorem}

Note that the choice of $\varphi$ is not unique, but one has the freedom of modifying it by harmonic terms -- a fact that reflects the gauge invariance of Schr\"odinger operators. In complex analysis, $\varphi$ plays the role of a weight function, thus $\varphi\ge0$ and $\triangle \varphi\ge0$ are reasonable assumptions. \\
In particular we have the following Corollary.

\begin{corollary}
\label{cor}
Suppose that the magnetic field is of the form $B=\triangle\varphi$ for some $\mathcal C^2$-function $\varphi$. If on the complement of a compact set it holds that 
$$B(x,y)\ge(x^2+y^2)^{-2} \quad \text{or}\quad B(x,y)\le-(x^2+y^2)^{-2},$$
then the resolvent of the Dirac operator is not compact .
\end{corollary}

In contrast to \cite{hnw}, our condition does not make assumptions on the derivatives of the magnetic field, but on its structure and growth. One can easily find examples of functions satisfying the assumptions of Theorem \ref{dirac} but not condition \eqref{helffer} and vice versa.\\
For instance, if $\varphi=x^3$, then $B$ is polynomial and thus satisfies  \eqref{helffer}, but $\varphi+h$ is never bounded from above or below for any harmonic function $h$. On the other hand, $\varphi=e^{x^2+y^2}$ satisfies the assumptions of Theorem \ref{dirac}, but each derivative of order $r+1$ of $\varphi$ grows faster at infinity than all derivatives of order $\le r$, hence \eqref{helffer} can not hold.\\
In some sense the two conditions are complementary. By Corollary \ref{cor}, the assumption of Theorem \ref{dirac} is satisfied if $B(x,y)\ge(x^2+y^2)^{-2}$ or $ B(x,y)\le-(x^2+y^2)^{-2}$, whereas the case $|B(x,y)|\to 0$ for $|(x,y)|\to\infty$ is covered by the result in \cite{hnw}.
\vskip 0.5cm

Our proof of Theorem \ref{dirac} relies on the following Theorem, which gives a characterization of discreteness of the spectrum for a certain class of differential operators. In particular it extends the main result in \cite{iwatsuka}. A variant of condition (3) in Theorem \ref{compact} recently appeared in \cite{has}. This was fitted to the $\dquer$-Neumann problem, nevertheless the ideas in the proof of \cite{has} actually led to this much more general statement.

\begin{theorem}
\label{compact}
Let $T$ be an invertible linear partial differental operator acting on $dom(T)$, which is closed, densely defined and elliptic in the interior of a smooth domain $\Omega\subseteq \Bbb R^n$. By this last property we mean that G\aa rding's inequality holds on each set relatively compact in $\Omega$. Let $T_\varphi^*$ be the adjoint of $T$ in $L^2(\Omega, \varphi)$ and set $P=T^*_\varphi T$.\\
Then the follwing are equivalent:
\begin{enumerate}
\item $P$ has compact resolvent, i.e., $(P-\lambda)^{-1}$ is compact for some (equivalently for all) $\lambda \in\Bbb C\setminus\sigma(P)$.
\item The injection $j_\varphi$ of $dom(T)$ equipped with the graph norm $u\mapsto\|Tu\|_\varphi$ into $L^2(\Omega, \varphi)$ is compact.
\item For all $\varepsilon>0$ there is $\Omega_\varepsilon\subset\subset\Omega$ such that $\| u\|_{L^2(\Omega\setminus\Omega_\varepsilon,\varphi)}<\varepsilon$ for all $u\in\mathcal L=\{ u\in dom(T):\ \|u\|_T<1\}$.
\item There is a smooth function $\lambda$, such that $\lambda\to\infty$ for $z\to\partial\Omega$ and
$$\langle Pu,u\rangle_\varphi\ge\int _\Omega\lambda| u|^2e^{-2\varphi}d\mu$$
for all $u\in dom(P)$, where $\mu$ denotes the Lebesgue measure on $\Omega$.
\end{enumerate}
\end{theorem}

{\bf Remark. } If $T$ is not invertible, one can look at $\ker(T)^\perp\subset L^2(\Omega, \varphi)$. In that case, one furthermore has to assume that $\dim\ker(T)<\infty$ in order to have compact resolvent for $P$.
\vskip 0.5 cm

{\bf Remark. } Theorem \ref{compact} generalizes the Main Theorem in \cite{iwatsuka}, where the same result was proven for magnetic Schr\"odinger operators with electric potentials that are semibounded from below.\\
If $\mathcal C_0^\infty(\Omega)$ is a core in the form domain (i.e., dense in the graph norm), one can push the analogy to \cite{iwatsuka} even further by also adding the bottom of the spectrum of the Dirichlet realization to the picture. This is, $P$ has compact resolvent if and only if the lowest eigenvalue of the Dirichlet realization of $P$ in $\Omega\setminus K_j$ tends to infinity as $j\to\infty$, for any sequence of smoothly bounded compact sets $K_j$ exhausting $\Omega$. Note that if $\Omega=\Bbb R^n$, then $\mathcal C_0^\infty$ is always a core. 
\vskip 0.5cm

From our proofs of the two previous Theorems, we get the following result as a Corollary. We point this out separately, as compact injection Theorems of this kind are of interest in statistics, as it turns out that they are essential in proving the existence of an orthonormal set of Nonlinear Principal Components.

\begin{corollary}
\label{sob}
Let $\varphi$ be a measurable weight function which is bounded from below. Set
$$H^1_\varphi(\mathbb{R}^n)=\{f\in  L^2(\mathbb R^n , \varphi)\ : \frac{\partial f}{\partial x_j}\in L^2(\mathbb R^n , \varphi)\ \forall\ 1\le j\le n\}$$
with the norm 
$$ \| f\|_{1,\varphi}^2=\| f\|_\varphi^2+\sum_{j=1}^n\left\|\frac{\partial f}{\partial x_j}\right\|^2_\varphi.$$
Then the injection $H^1_{\varphi }(\mathbb R^n)\hookrightarrow L^2(\mathbb{R}^n, \varphi)$ is never compact.
\end{corollary}
\vskip 1cm

\section{Proof of Theorem \ref{compact}.}~\\

In order to prove the Theorem, we will make use of the following functional analytic characterization of precompact sets in weighted Lebesgue spaces. 

\begin{lemma}
\label{precp}
Let $\Omega\subset \Bbb R^n$ and $\mathcal A$ be a bounded subset of $L^2 (\Omega , \varphi)$. Then $\mathcal A$ is precompact if and only if the following two conditions hold:
\end{lemma}
\begin{enumerate} 
\item
for all $\varepsilon  > 0$ and all $\Omega^\prime\subset\subset\Omega$ there exists $\delta > 0$ such that 
$$\| \tau_h f-f\|_{L^2(\Omega^\prime,\varphi)}<\varepsilon $$
for each $h \in \Bbb R^n$ with $|h| <\delta$ and all $f \in\mathcal A$, where $\tau _h f (x ) = f (x + h)$.
\item 
for all $\varepsilon > 0$ there exists $\Omega_\varepsilon\subset\subset\Omega$ such that 
$$\| f\|_{L^2(\Omega\setminus\Omega_\varepsilon ,\varphi)}<\varepsilon$$
for each $f\in\mathcal  A$. 
\end{enumerate}

For the proof we refer to \cite{adams}, Theorem 2.32. See also \cite{bre}, Corollaire IV.26.\\

{\em Proof of Theorem \ref{compact}.} Let $P^{-1}$ be the inverse of $P$ and let $j_\varphi$ be the injection of $dom(T)$ into $L^2(\Omega, \varphi)$. Equip $dom(T)$ with the graph norm by setting the inner product to be $\langle u,v\rangle_T=\langle Tu,Tv\rangle_\varphi$. We first show that $P^{-1}=j_\varphi\circ j_\varphi^*$, inspired by an idea of E. Straube in \cite{str}. For all $u,v\in dom(T)$ it holds
$$\langle u,v\rangle_\varphi=\langle u,j_\varphi v\rangle_\varphi=\langle j_\varphi^* u,v\rangle_T,$$
while on the other hand
$$\langle u,v\rangle_\varphi=\langle PP^{-1}u,v\rangle_\varphi=\langle TP^{-1}u,Tv\rangle_\varphi=\langle P^{-1}u,v\rangle_T.$$
Hence, $P^{-1}=j_\varphi^*$ as an operator to $dom(T)$ and consequently  $P^{-1}=j_\varphi\circ j_\varphi^*$ as an operator to $L^2(\Omega, \varphi)$. This proves the equivalence of $(1)$ and $(2)$.\\

Now we show that $(2)\implies (3)\implies(4)\implies (2)$.\\ 
Suppose that the injection $j_\varphi$ is compact. Hence the image of $\mathcal L$ is precompact in $L^2(\Omega, \varphi)$, thus $(2)\implies (3)$ by Lemma \ref{precp}.\\

If $(3)$ holds, then by linearity of $T$ for all $\varepsilon>0$ there is $\Omega_\varepsilon\subset\subset\Omega$, such that $\| u\|_{L^2(\Omega\setminus\Omega_\varepsilon)}\le\varepsilon\| u\|_T$ for all $u\in dom(T)$. Thus for all $u\in dom(P)\subset dom(T)$:
\begin{align*}
\int\limits_\Omega | u|^2 e^{-2\varphi}d\mu\le &\int \limits_{\Omega\setminus\Omega_\frac{1}{4}}1\cdot | u|^2 e^{-2\varphi}d\mu+\int \limits_{\Omega_\frac{1}{4}\setminus\Omega_\frac{1}{16}} 2\cdot | u|^2 e^{-2\varphi}d\mu+\int \limits_{\Omega_\frac{1}{16}\setminus\Omega_\frac{1}{64}}4\cdot | u| ^2e^{-2\varphi}d\mu+\dots\\
\le&\ 2\ \|u\|_T^2.
\end{align*}
Hence it is clear that one can find a smooth function $\lambda$ tending to infinity at the boundary of $\Omega$ such that 
$$\langle Pu,u\rangle_\varphi=\|u\|_T^2 \ge\int_\Omega \lambda| u|^2e^{-2\varphi}d\mu$$
for all $u\in dom(P)$.\\

Finally suppose that $(4)$ holds. We will prove $(2)$ by checking the two conditions from Lemma \ref{precp} for the unit ball $\mathcal L$ in $dom(T)$. Since by smoothness of $ \Omega$ functions in $\mathcal C^\infty(\overline\Omega)$ are dense in $dom(T)$ in the graph norm, we can restrict ourselves to these. Now if $\omega \subset\subset \Omega,$ choose $\omega \subset\subset\omega_1 \subset\subset \omega_2 \subset\subset \Omega $ and a smooth cut-off function $\psi $ with $\psi (z)=1$ for $z\in \omega_1$ and $\psi (z)=0$ for $z\in \Omega \setminus \omega_2 .$ For $u\in dom(T)$, define $\tilde u = \psi u$ and note that the domain of $T$ is preserved under multiplication by a function in $\mathcal{C}^\infty_0(\Omega )$. Therefore $\tilde u$ has compactly supported coefficients and belongs to $dom(T).$ The graph norm of $\tilde u$ is bounded by a constant $C$ depending only on $\omega, \omega_1, \omega_2, \psi, $ if $u$ belongs to $\mathcal L$. By construction we have
 $$\| \tau_h u -u \|_{ L^2(\omega,\varphi )}=\| \tau_h \tilde u -\tilde u \|_{ L^2(\omega,\varphi )}$$
for $|h|<dist(\omega,\omega_1)$, hence it sufficies to estimate the latter expression $\| \tau_h \tilde u -\tilde u \|_{ L^2(\omega,\varphi )}$. Since $T$ is elliptic in the interior, this essentially comes from G\aa rding's inequality which says that for any smoothly bounded domain $\Omega^\prime\subset\subset\Omega$ there is $C_{\Omega^\prime,\varphi}>0$ such that 
$$\| u\|_{1,\varphi}^2\le  C_{\Omega^\prime,\varphi}\left(\|Tu\|^2_\varphi+\|u\|^2_\varphi\right)$$
for all $u\in\mathcal C_0^\infty(\Omega^\prime)$. So let $\tilde u\in\mathcal C_0^\infty(\omega_2)$ and set for $h\in\Bbb R^n$ and $t\in\Bbb R$
$$ \tilde v(t)=\tilde u(x+ht).$$
Note that 
$$|\tilde v'(t)|\le |h| \left [ \sum_{k=1}^n \left | \frac{\partial\tilde u}{\partial x_k}(x+th)\right |^2 \right ]^{1/2}.$$
By the fact that 
$$\tilde u(x+h)-\tilde u(x)=\tilde v(1)-\tilde v(0)=\int_0^1\tilde v_j'(t)\,dt,  $$
 we can estimate for $|h|<1/2\ dist(\omega,\omega_1)$
 \begin{align*}
 &\int_{\omega}|\tau_h\tilde u(x)-\tilde u(x)|^2 e^{-2\varphi (x)}\,d\mu (x)\\
  \le& |h|^2 \,  \int_{\omega}\left [ \int_0^1 \sum_{k=1}^n\left | \frac{\partial\tilde u}{\partial x_k}(x+th)\right |^2\,dt \, \right ] e^{-2\varphi (x)} \,d\mu (x) \\
\le& C_{\omega_1,\varphi}\,  |h|^2 \,  \int_{\omega_2} \sum_{k=1}^n \left | \frac{\partial\tilde u}{\partial x_k}(x)\right |^2 e^{-2\varphi (x)} \,d\mu (x) .
\end{align*}
Now, since the cut-off function $\psi$ was fixed, we get by using G\aa rding's inequality
\begin{align*}
\Vert \psi u \Vert ^2_{1, \varphi}&\leq C'_{\varphi,\omega_2}
\left(\Vert T ( \psi u)\Vert^2_\varphi+\Vert \psi u\Vert^2_\varphi \right) \\
&\leq C''_{\varphi,\omega_2,\psi} \left(\Vert T u\Vert^2_\varphi+\Vert u\Vert^2_\varphi \right) .
\end{align*}
By assumption, $T^{-1}$ is bounded, which is equivalent to $\| u\|_\varphi\lesssim\|u\|_T$. So we can neglect the second term and summing up, we get
$$\| \tau_h u -u \|_{ L^2(\omega,\varphi )}\le C_{\varphi, \omega_2, \psi}'''|h|^2\|Tu\|_\varphi^2.$$
Since we started from the unit ball $\mathcal L$ in $dom(T)$, we get that condition (i) of Lemma \ref{precp} is satisfied.\\
Now we verify condition (ii). Let $\varepsilon>0$ be given and choose $M$ such that $1/M\le\varepsilon$ and $\Omega_M\subset\subset\Omega$ such that $\lambda\ge M$ on $\Omega\setminus\Omega_M$.  Thus we have for all $u\in\mathcal L\cap dom(P)$
\begin{align*}
\|u\|_T\ge\int_\Omega\lambda| u|^2e^{-2\varphi}d\mu\ge M\int_{\Omega\setminus\Omega_M}| u|^2e^{-2\varphi}d\mu,
\end{align*}
which makes (ii) immediate, since $dom(P)$ is dense in $dom(T)$.
\begin{flushright}
$\square$
\end{flushright}

{\bf Remark. } We chose that way of proving the implication $(4)\implies (2)$ since we think it shows most clearly what is going on. A shorter way using more functional analysis, in particular the Rellich -- Kondrachov Theorem, is the following: An operator $K:H_1\to H_2$ between two Hilbert spaces is compact if and only if for each $\varepsilon>0$ there is a compact operator $K_\varepsilon:H_1\to H_2$ such that 
$$\| Kf\|_{H_2}^2\le\varepsilon \|f\|_{H_1}^2+C_\varepsilon\|K_\varepsilon f\|_{H_2}^2.$$
Hence, starting from $(4)$ we immediately get
\begin{align*}
\| j_\varphi f\|_\varphi^2=\int_\Omega |f|^2e^{-2\varphi}d\mu\le&\frac{1}{M}\int_{\Omega\setminus\Omega_M}\lambda |f|^2e^{-2\varphi}d\mu+C_M\int_{\Omega_M}|f|^2e^{-2\varphi}d\mu\\
\le&\frac{1}{M}\|f\|^2_T+C^\prime_M\|f\|^2_{L^2(\Omega_M)}.
\end{align*}
Now the injection $dom(T)\hookrightarrow L^2(\Omega_M)$ is compact by G\aa rding's inequality combined with the Rellich -- Kondrachov Theorem, which shows $(2)$.
\vskip 1cm

\section{Proof of Theorem \ref{dirac}.}~\\

As a first step, we use the following reduction of the problem to Pauli operators, see \cite{hahe}.  Suppose that $\Bbb D $ has compact resolvent. Then also $\Bbb D^2$ has, since 
$$(\Bbb D^2-i)^{-1}=(\Bbb D+\sqrt{i})^{-1}\circ(\Bbb D-\sqrt i)^{-1}.$$
Note that by a classical result, the spectrum of $\Bbb D$ is disjoint to $(-1,1)$, see e.g. \cite{hnw} for a proof. So $\Bbb D$ is always invertible.
An easy computation shows the standard fact that for the square of the Dirac operator it holds
\begin{equation*}
\Bbb D^2=
\begin{pmatrix}
P_+ &0\\
0 &P_-\\
\end{pmatrix},
\end{equation*}
where $P_\pm$ are the so-called Pauli operators,
\begin{equation}
P_\pm =- \left(\frac{\partial}{\partial x}-iA_1(x,y)\right)^2-\left(\frac{\partial}{\partial y}-iA_2(x,y)\right)^2\pm B(x,y).
\end{equation}
This implies that if $\Bbb D$ has compact resolvent, then both $P_\pm$ have (Now as operators on $L^2(\Bbb R^2,\varphi)$). We shall show that if $e^{-2\varphi}$ is bounded, the resolvent of $P_-$ is never compact. In the case that $e^{2\varphi}$ is bounded, one needs to replace $\varphi$ by $-\varphi$ and notice that $P_+$ and $P_-$ swap their roles.
\vskip 0.5cm

Now consider the weighted space 
$$L^2(\mathbb C , \varphi )=\{ f:\mathbb C \to \mathbb C \ : \ \int_{\mathbb C}
|f|^2\, e^{-2\varphi}\,d\mu < \infty \},$$
where $\mu$ denotes the Lebesgue measure on $\Bbb C$. It is well-known in complex analysis, that 
$$\dquer=\frac{\partial}{\partial \overline z}=\frac{1}{2}\left(\frac{\partial}{\partial x}+i\frac{\partial}{\partial y}\right)$$
is a locally elliptic closed and densely defined operator on $L^2(\mathbb C , \varphi )$. Identifying $\Bbb C\simeq\Bbb R^2$ and defining an operator $\overline D $ on $L^2(\Bbb R^2)$,
$$\overline D=e^{-\varphi}\dquer e^\varphi,$$
we see from the unitary equivalence, that $\dquers\dquer$ has compact resolvent if and only if $\overline D^*\overline D$ has. Since $4\overline D^*\overline D=P_-$ as was pointed out in \cite{hahe}, the proof of Theorem \ref{dirac} boils down to showing that the injection $dom(\dquer)\hookrightarrow L^2(\mathbb C , \varphi )$ is not compact and using Theorem \ref{compact}.
\vskip 0.5cm

So suppose that $P_- $ has compact resolvent. By the above consideration and Theorem \ref{compact}, this is the case if and only if there is some smooth function $\lambda$ with $\lambda\to\infty$ as $|z|\to\infty$ such that 
$$\int_{\Bbb C}\left|\frac{\partial f}{\partial \overline z}\right|^2e^{-2\varphi}d\mu\ge\int_{\Bbb C} \lambda| f|^2e^{-2\varphi}d\mu$$
for all $f\in\mathcal C_0^\infty(\Bbb C)$. Here we used that we can restrict ourselves to $\mathcal C_0^\infty(\Bbb C)$ since it is dense in $dom(\dquer)\cap dom(\dquers)$ in the graph norm, see \cite{gh} Lemma 2.2. By Definition, $\dquers \dquer $ is a positive operator, hence we can assume $\lambda\ge0$. Without loss of generality we also assume that $\lambda\ge\varepsilon>0$ in a neighborhood $U$ of $0$.\\
Let $\{\chi_R\}_{R\in\Bbb N}$ be a family of smooth cut-off functions which are identically one on $\Bbb B_R$, the ball with radius $R$ and center $0$, supported in $\Bbb B_{R+1}$  and have uniformly bounded first order derivatives. In fact, we can assume that $\sup|\nabla\chi_R(z)|\le2$ for all $z\in\Bbb C$ and $R\in \Bbb N$. Now $\chi_R \in \mathcal C_0^\infty(\Bbb C)\subset dom(\dquers\dquer)$, thus we have
$$\int\limits_{supp( \nabla\chi_R)}\left|\frac{\partial \chi_R}{\partial \overline z}\right|^2e^{-2\varphi}d\mu\ge\int\limits_{supp(\chi_R)} \lambda| \chi_R|^2e^{-2\varphi}d\mu.$$
Using the assumption on the derivatives of $\chi_R$, we get
$$4\int\limits_{\Bbb B_{R+1}\setminus\Bbb B_R}e^{-2\varphi}d\mu\ge\int\limits_{\Bbb B_R} \lambda |\chi_R|^2e^{-2\varphi}d\mu\ge\varepsilon\int\limits_{U}e^{-2\varphi}d\mu>\delta>0,$$
yielding that 
$$\| 1\|^2_\varphi=\sum_{R\in \Bbb N}\ \int\limits_{\Bbb B_{R+1}\setminus\Bbb B_R}e^{-2\varphi}d\mu$$
can not be finite.\\
On the other hand, if the injection $dom(\dquer)\hookrightarrow L^2(\mathbb C , \varphi )$ is compact, then the same holds true for $dom(\nabla)\hookrightarrow L^2(\mathbb C , \varphi )$, simply by triangle inequality. Note that $dom(\nabla)=H^1_\varphi(\Bbb C),$ where 
$$H^1_\varphi(\Bbb C)=\{f\in  L^2(\mathbb C , \varphi)\ : \frac{\partial f}{\partial x},  \frac{\partial f}{\partial y}\in L^2(\mathbb C , \varphi)\}.$$
Now the proof is finished by citing Theorem 3.3 from \cite{an}, where it was shown by using a method from \cite{adams} that 
$$\int_{\Bbb C}e^{-2\varphi}d\mu<\infty$$
is a necessary condition for compactness of the injection $H^1_\varphi(\Bbb C)\hookrightarrow L^2(\mathbb C , \varphi )$, when the weight function is bounded from above.
\begin{flushright}
$\square$
\end{flushright}

{\em Proof of Corollary \ref{cor}. } We only need to consider the case $B(x,y)\ge(x^2+y^2)^{-2}$, else $B\mapsto -B$.\\
Since $\dquer$ annihilates holomorphic functions, we can replace $\chi_R$ in the proof of Theorem \ref{dirac} by $\chi_R h$ for any entire function $h$ and see by the same arguments, that if $P_-$ has compact resolvent, there can not be a holomorphic function which is integrable with respect to the weight $\varphi$.\\
On the other hand, the condition $\triangle\varphi\ge(x^2+y^2)^{-2}$ assures that $\dim\mathcal H(\Bbb C)\cap L^2(\Bbb C,\varphi)=\infty$ (see Lemma 3.4 in \cite{shi}), a contradiction.
\begin{flushright}
$\square$
\end{flushright}

{\em Proof of Corollary \ref{sob}. } Let $T=d$ be defined as 
$$Tf=\sum_{j=1}^n\frac{\partial f}{\partial x_j}dx_j$$
and equip the space of $1$-forms with the inner product equal the sum of the coefficient-wise weighted inner products. Then $dom(T)=H^1_\varphi(\Bbb R^n)$ and by Definition, G\aa rding's inequality holds for $T$. Hence by Theorem \ref{compact}, the injection $H^1_\varphi(\Bbb R^n)\hookrightarrow L^2(\Bbb R^n , \varphi )$ is compact if and only if 
$$\sum_{j=1}^n\int_{\Bbb R^n}\left|\frac{\partial f}{\partial x_j}\right|^2e^{-2\varphi}d\mu\ge\int_{\Bbb R^n} \lambda| f|^2e^{-2\varphi}d\mu$$
for a function $\lambda$ with $\lambda\to\infty$ for $|z|\to\infty $ and all $f \in \mathcal C^\infty_0(\Bbb R^n)$. By the arguments given in the proof of Theorem \ref{dirac}, this implies $1\notin L^2(\mathbb R^n , \varphi )$. But if the weight function is bounded, $1\in L^2(\mathbb R^n , \varphi )$ is a necessary condition for compacteness of this embedding by Theorem 3.3 in \cite{an}, a contradiction.
\begin{flushright}
$\square$
\end{flushright}

{\bf Aknowledgement.} The author is thankful to B. Helffer for useful and interesting discussions on the Dirac operator. \\
Moreover he wants to thank F. Haslinger for discussions on the complex analytic part of this work and in particular for drawing his attention at Lemma \ref{precp} and its connection to compactness. These conversations actually led to the ideas used in the proof of Theorem \ref{compact}.

\end{document}